\documentclass[reqno,12pt]{amsart}

\usepackage{geometry}
\usepackage{epsfig}
\usepackage{amsmath}
\usepackage{amsfonts}
\usepackage{mathrsfs}
\usepackage{cite}
\numberwithin{equation}{section}
\allowdisplaybreaks[4]

\begin{document}
\newtheorem{theorem}{Theorem}[section]
\newtheorem{lemma}[theorem]{Lemma}
\newtheorem{corollary}[theorem]{Corollary}
\newtheorem{conjecture}[theorem]{Conjecture}
\newtheorem{qn}[theorem]{Question}
\newtheorem{pro}[theorem]{Proposition}
\newtheorem{remark}[theorem]{Remark}
\newcommand{\arcsinh}{\mbox{arcsinh}}
\newcommand{\arccosh}{\mbox{arccosh}}
\renewcommand{\theequation}{\arabic{section}.\arabic{equation}}
\newcommand{\IN}{\mathbb{N}}
\newcommand{\IC}{\mathbb{C}}
\newcommand{\e}{\mbox{e}}
\newcommand{\arth}{\mbox{arth}}
\newcommand{\K}{\mathscr{K}}
\newcommand{\E}{\mathscr{E}}
\newcommand{\s}{\sum_{n=0}^{\infty}}
\newcommand{\su}{\sum_{n=1}^{\infty}}
\newcommand{\la}{\lambda}
\newcommand{\G}{\Gamma}
\newtheorem{question}[theorem]{Question}

\def\square{\hfill${\vcenter{\vbox{\hrule height.4pt \hbox{\vrule width.4pt
height7pt \kern7pt \vrule width.4pt} \hrule height.4pt}}}$}

\title[]
{Sharp double inequality for complete\\ elliptic integral of the first kind}


\author{Qi Bao$^{\ast}$}

\address{Qi Bao, School of Mathematical Sciences, East China Normal University, Shanghai, 200241, China}

\email{52205500010@stu.ecnu.edu.cn}



\thanks{$^{\ast}$Corresponding author.}

\maketitle
\begin{abstract}
For $r\in(0,1)$, the function $\K(r)=\int_0^{\pi/2}(1-r^2\sin^2t)^{-1/2}dt$ is known as the complete elliptic integral of the first kind. In this paper, we prove the absolute monotonicity of  two functions involving $\K(r)$. As a consequence, we improve Alzer and Richards' result.
\end{abstract}

\bigskip
\section{introduction}
In the past few centuries, the complete elliptic integral of the first
kind (cf. \cite{AS,Anderson1,QVu1}) $\K(r)$ defined on $(0,1)$ by
\begin{align}
\K=\K(r)=\int_0^{\pi/2} \frac{1}{\sqrt{1-r^2\sin^2 t}} dt
=\frac{\pi}{2}F\left(\frac12,\frac12;1;r^2\right),
\end{align}
where $F$ denotes the classical Gaussian hypergeometric function (cf. \cite{AQVV,Carlson})
\begin{align}\label{Gauss}
F(a,b;c;x)={}_2F_1(a,b;c;x)=\s\frac{(a,n)(b,n)}{(c,n)n!}x^n, \quad |x|<1,
\end{align}
where $(a,n)$ is the Pochhammer symbol or shifted factorial defined as $(a,0)=1$ for $a\neq0$, and
\begin{align}
(a,n)=a(a+1)(a+2)\cdots(a+n-1)=\frac{\G(n+a)}{\G(a)}
\end{align}
for $n\in\IN=\{1,2,3,\cdots\}$, where
\begin{align}
\G(x)=\int_0^{\infty}t^{x-1}e^{-t}dt \quad (x>0)
\end{align}
is the classical Euler Gamma function \,(cf. \cite{AS}).

It is well known that the complete elliptic integrals have many important
applications in physics, engineering, geometric function theory, quasiconformal analysis, theory of mean values, number theory and other related fields\, (cf
\cite{AbbasBaloch,Ramanujan,BSF,QAV1,AVV2,ByrdFriedman,WGQ2,AlzerQiu1,BRP,HTC1,SHS1,AQV5,AVV6,BF2}).

Recently, the complete elliptic integrals have attracted the attention of numerous mathematicians. It is well known that complete elliptic integrals cannot be represented by the elementary transcendental functions. Therefore, there is a need for sharp computable bounds for the family of integrals. In particular, many remarkable properties and inequalities for the complete elliptic integrals can be found in the literature \cite{QAV1,WCLC1,AV3,AlzerQiu1,AlzerRichards,GQ1,HTM2,YTZ2}.

For example, in order to refine the following well-known asymptotic formula
\begin{align}\label{asymptotic}
\lim_{r\to1^-}\left[\K(r)-\log \left(\frac{4}{\sqrt{1-r^2}}\right)\right]=0,
\end{align}
Anderson, Vamanamurthy and Vuorinen in \cite{AV3} conjectured that the inequality
\begin{align}\label{In5}
\K(r)<\log \left(\frac{4}{\sqrt{1-r^2}}\right)
-\left(\log5-\frac{\pi}{2}\right)(1-r)
\end{align}
holds for each $r\in(0,1)$. Later, the conjecture was proved by Qiu et al. in \cite{QAV1}.

In 2020, Wang, Chu, Li and Chu in \cite{WCLC1} improved (\ref{In5}), and showed that the inequality
\begin{align}\label{In6}
&\log \left( 1+\frac{4}{\sqrt{1-r^2}}\right)- \left(\log5-\frac{\pi}{2}\right)
+\left(\frac{\pi}{8}-\frac25\right)r^2+\alpha^* r^4<\K(r) \nonumber \\
&<\log \left( 1+\frac{4}{\sqrt{1-r^2}}\right)- \left(\log5-\frac{\pi}{2}\right)
+\left(\frac{\pi}{8}-\frac25\right)r^2+\beta^* r^4
\end{align}
is valid for $r\in(0,1)$ with the best possible constants $\alpha^*=9\pi/128-11/50=0.000893\cdots$ and $\beta^*=2/5+\log5-5\pi/8=0.0459\cdots$.

Very recently, Alzer and Richards in \cite{AlzerRichards} obtained the following upper bound for $\K(r)$, that is,
\begin{align}\label{In2}
\K(r)<\frac{\pi}{2} \frac{16-5\log(1-r^2)}{16+(5\pi-16)r^2} \,\, (0<r<1).
\end{align}
Observe that the upper bound of $\K(r)$ in (\ref{In2}) is concise and meaningful. Based on the known results such as those above mentioned, the following question is natural:
\begin{question}\label{Question}
 Whether we can find an improved upper bound  and  a similar form of  lower bound for (\ref{In2})?
\end{question}

The main purpose of this paper is to give a positive answer to Question \ref{Question}. Our results are following Theorems \ref{theorem1}-\ref{theorem2} and Corollary \ref{corollary}. For convenience, we let $\IN$ denotes the set of positive integers as usual, put
\begin{align}
\theta&=\frac{\pi(17-5\pi)}{32}=0.126845\cdots,\label{theta}\\
\lambda&=\frac85-\log4=0.213705\cdots,\label{lambda}\\
\alpha&=\frac{85}{8}\pi-\frac{185}{32}\pi^2+\frac{25}{32}\pi^3
=0.544425\cdots,\label{alpha}\\
\beta&=(8-10\log2)\pi-\frac{85}{32}\pi^2+\frac{25}{32}\pi^3
=1.364397\cdots,\label{beta}\\
\delta&=\frac{128}{5}-32\log2-\frac{17}{2}\pi+\frac{5}{2}\pi^2
=1.389763\cdots,\label{delta}\\
\zeta&=-\frac{128}{5}+32\log2+\left(\frac{47}{8}-10\log2\right)\pi
  +\frac{5}{8}\pi^2=-0.569791\cdots. \label{zeta}
\end{align}

\begin{theorem}\label{theorem1}
Let $\theta$, $\alpha$ and $\beta$ are given in
(\ref{theta}), (\ref{alpha}) and (\ref{beta}), respectively. Define the function $f$ on $(0,1)$ by
\begin{align*}
f(r)=\frac{\pi}{2} \left[16-5\log(1-r^2)\right]
-\left[\theta r^2+\K(r)\right]\left[16+(5\pi-16)r^2\right].
\end{align*}
Then all coefficients are positive in the Maclaurin series for $f_1 \equiv f/r^4$ in powers of $r^2$ with range $(\alpha,\beta)$. In other words, $f_1$ is absolutely monotonic on $(0,1)$.

\end{theorem}

\begin{theorem}\label{theorem2}
Let $\lambda$ and $\delta$ are given in
(\ref{lambda}) and (\ref{delta}), respectively. Define the function $g$ on $(0,1)$ by
\begin{align*}
g(r)=\left[\lambda r^2+\K(r)\right]
\left[16+(5\pi-16)r^2\right]-\frac{\pi}{2} \left[16-5\log(1-r^2)\right].
\end{align*}
Then all coefficients are negative in the Maclaurin series for $g_1 \equiv g/r^2$ in powers of $r^2$ with range $(0,\delta)$. In other words, $-g_1$ is absolutely monotonic on $(0,1)$.
\end{theorem}

\begin{corollary}\label{corollary}\rm{}
According to Theorem \ref{theorem1}-\ref{theorem2}, we can find better bounds for $\K(r)$ than (\ref{In2}). For example, Theorem \ref{theorem1} (Theorem \ref{theorem2}) implies that the function $f_1$\,($g_1$) is is strictly increasing and convex (decreasing and concave) from $(0,1)$ onto $(\alpha,\beta)$\, ($(0,\delta)$, respectively). Consequently, the double inequality
\begin{align}\label{In7}
\max \Bigg\{& \frac{\pi[16-5\log(1-r^2)]-2[\alpha+(\beta-\alpha)r] r^4}{32+2(5\pi-16)r^2}-\theta r^2, \nonumber \\
&\frac{\pi[16-5\log(1-r^2)]+2\delta (1-r)r^2}{32+2(5\pi-16)r^2}-\lambda r^2 \Bigg\}
\leq \K(r)\leq \\
&\min \Bigg\{\frac{\pi[16-5\log(1-r^2)]-2\alpha r^4}{32+2(5\pi-16)r^2}-\theta r^2,
\, \frac{\pi[16-5\log(1-r^2)]+2\delta r^2}{32+2(5\pi-16)r^2}-\lambda r^2\Bigg\}.
\nonumber
\end{align}
holds for all $r\in(0,1)$.  The first (second) equality holds if and only if $r\to 0$ or $r\to1$ ($r\to 0$, respectively).
\end{corollary}

\section{proof of theorem \ref{theorem1}-\ref{theorem2}}
In order to prove our main results in this section, we need next lemma.
\begin{lemma}\label{lemma}
For $n\in \IN$, the sequence
\begin{align*}
Q_n=\frac{5\pi n}{5\pi n^2-16n+4}-\left[\frac{\G(n-1/2)}{\G(n)}\right]^2
\end{align*}
is positive.
\end{lemma}

\proof
In \cite[equation (1.3)]{Kershaw1}, Kershaw proved that
\begin{align}\label{In3}
\left(x+\frac{s}{2}\right)^{1-s}<\frac{\G(x+1)}{\G(x+s)}<
\left[x-\frac12+\left(\frac14+s\right)^{1/2}\right]^{1-s}
\end{align}
holds for $x>0$ and $0<s<1$. Hence by  first inequality sign of (\ref{In3}) and by using the substitution $x=n-1$, take $s=1/2$, we have
\begin{align}\label{In10}
\left[\frac{\G(n-1/2)}{\G(n)}\right]^2<\frac{4}{4n-3}
\end{align}
is vaild for $n\in \IN$. It follows from (\ref{In10}) that
\begin{align}\label{QnPn}
Q_n> \frac{5\pi n}{5\pi n^2-16n+4}-\frac{4}{4n-3}=P_n.
\end{align}
It is enough to prove  $P_n>0$ for $n\in \IN$, as a matter of fact,
\begin{align*}
P_n>0 &\Longleftrightarrow 5\pi n(4n-3)-4(5\pi n^2-16n+4)>0 \\
&\Longleftrightarrow (64-15\pi)n-16> 0
\end{align*}
holds for $n\in \IN$. Therefore, together with (\ref{QnPn}), yields the sequence $\{Q_n\}$ is positive for $n\in \IN$.
\endproof

{\it Proof of Theorem \ref{theorem1}}

By (\ref{Gauss}), expanding in power series yields
\begin{align*}
f(r)
&=\frac{\pi}{2} \left(16+5\su \frac1n r^{2n}\right)
-[16+(5\pi-16)r^2]
\left( \theta r^2+\frac{\pi}{2} \s \frac{(1/2,n)^2}{(n!)^2}r^{2n} \right) \\
&=8\pi-16\theta r^2-\theta(5\pi-16)r^4 \\
&\quad +\frac{5\pi}{2} \su \frac1n r^{2n}
-\frac{(5\pi-16)\pi}{2} \su \frac{(1/2,n-1)^2}{[(n-1)!]^2}r^{2n}
-8\pi \s \frac{(1/2,n)^2}{(n!)^2}r^{2n}\\
&= \left(\frac{85}{8}\pi-\frac{185}{32}\pi^2+\frac{25}{32}\pi^3\right)r^4\\
&\quad +\sum_{n=3}^{\infty} \left(
\frac{5\pi}{2n}-\frac{(5\pi n^2-16n+4)\G(n-1/2)^2}{2\G(n+1)^2}\right) r^{2n}\\
&= \alpha r^4+\sum_{n=3}^{\infty} \frac{(5\pi n^2-16n+4)Q_n}{2n^2}r^{2n},
\end{align*}
where $Q_n$ is given in lemma \ref{lemma}. Hence, $f_1=f/r^4$ has following series expansion
\begin{align}\label{f1}
f_1(r)=\alpha+\sum_{n=3}^{\infty} \frac{(5\pi n^2-16n+4)Q_n}{2n^2}r^{2n-4}.
\end{align}
It is easy to verify that $5\pi n^2-16n+4>0$ for $n\in \IN$, hence it follows from Lemma \ref{lemma} that all coefficients are positive in the Taylor series for $f_1$ in powers of $r^2$. By (\ref{f1}), we clearly see that
\begin{align}\label{f1(0)}
f_1(0^+)=\lim_{r\to0}f_1(r)=\alpha.
\end{align}
From (\ref{asymptotic}) it is easy to obtain the following asymptotic formula
\begin{align}\label{asy}
\K(r)\sim \log4-\frac{\log(1-r^2)}{2}, \quad  {\rm as} \,\, r\rightarrow 1.
\end{align}
Hence we obtain
\begin{align}\label{f1(1)}
f(1^-)&=\lim_{r\to1}f_1(r) \nonumber \\
&=\lim_{r\to1} \Bigg\{ 8\pi-\theta \left[16+(5\pi-16)r^2\right]r^2 \nonumber\\
&\quad -\frac{5\pi}{2}\log(1-r^2)
-\left(\log4-\frac{\log(1-r^2)}{2}\right)\left[16+(5\pi-16)r^2\right] \Bigg\} \nonumber\\
&=(8-10\log2)\pi-\frac{85}{32}\pi^2+\frac{25}{32}\pi^3=\beta.
\end{align}
Therefore, Theorem \ref{theorem1} directly follows from (\ref{f1(0)}) and  (\ref{f1(1)}) together with (\ref{f1}). This completes the proof.

{\it Proof of Theorem \ref{theorem2}}

Similarly, it follows from (\ref{Gauss}) that
\begin{align*}
g(r)
&=\left[16+(5\pi-16)r^2\right]\left(\lambda r^2+\frac{\pi}{2} \s\frac{(1/2,n)^2}{(n!)^2}r^{2n}\right)
-\frac{\pi}{2} \left(16+5 \su \frac1n r^{2n}\right)\\
&= -8\pi+16\lambda r^2+\lambda (5\pi-16)r^4 \\
&\quad -\frac{5\pi}{2} \su \frac1n r^{2n}
+\frac{(5\pi-16)\pi}{2} \su \frac{(1/2,n-1)^2}{[(n-1)!]^2}r^{2n}
+8\pi \s \frac{(1/2,n)^2}{(n!)^2}r^{2n}\\
&= \left(\frac{128}{5}-32\log2-\frac{17}{2}\pi+\frac{5}{2}\pi^2\right)r^2 \\
&\quad + \left[-\frac{128}{5}+32\log2+\left(\frac{47}{8}-10\log2\right)\pi
  +\frac{5}{8}\pi^2\right]r^4\\
&\quad -\sum_{n=3}^{\infty} \left(
\frac{5\pi}{2n}-\frac{(5\pi n^2-16n+4)\G(n-1/2)^2}{2\G(n+1)^2}\right) r^{2n}\\
&=\delta r^2+\zeta r^4
-\sum_{n=3}^{\infty} \frac{(5\pi n^2-16n+4)Q_n}{2n^2}r^{2n},
\end{align*}
where $Q_n$ is given in Lemma \ref{lemma}. Hence, $g_1=g/r^2$ has following series expansion
\begin{align}\label{g1expand}
g_1(r)=\delta+\zeta r^2
-\sum_{n=3}^{\infty} \frac{(5\pi n^2-16n+4)Q_n}{2n^2}r^{2n-2}.
\end{align}
Therefore, it follows from Lemma \ref{lemma} that all coefficients are negative in the Taylor series for $g_1$ in powers of $r^2$. By (\ref{g1expand}), we clearly see that
\begin{align}\label{g1(0)}
g_1(0^+)=\lim_{r\to0} g_1(r)=\delta.
\end{align}
Again using asymptotic formula (\ref{asy}), we clearly see that
\begin{align}\label{g1(1)}
g_1(1^-)&=\lim_{r\to 1} g_1(r) \nonumber \\
&=\lim_{r\to 1}\Bigg\{ -8\pi+\lambda \left[16+(5\pi-16)r^2\right]r^2 \nonumber\\
&\quad +\frac{5\pi}{2}\log(1-r^2)
+\left[16+(5\pi-16)r^2\right] \left(\log4-\frac{\log(1-r^2)}{2}\right) \Bigg\}
=0.
\end{align}

Therefore, Theorem \ref{theorem2} directly follows from
(\ref{g1expand})-(\ref{g1(1)}). This completes the proof.
\square

\begin{remark}\rm{}
Clearly, the upper bound of inequality (\ref{In7}) is better than (\ref{In2}). Moreover, computer simulation and experiments of Maple 2016 show that (\ref{In6}) and (\ref{In7})  have their own merits.
\end{remark}

\end{document}